\documentclass[reqno]{amsart}
\usepackage{hyperref}
\usepackage{setspace}
\usepackage{graphicx}
\usepackage{amssymb}
\usepackage{mathrsfs}

\begin{document}
    \title[\hfilneg ]
    { EXISTENCE AND UNIQUENESS OF  NONLOCAL
        BOUNDARY CONDITIONS FOR
        HILFER-HADAMARD-TYPE FRACTIONAL
        DIFFERENTIAL EQUATIONS}
 \author[\hfil\hfilneg]{Ahmad Y. A. Salamooni, D. D. Pawar }
 \address{Ahmad Y. A. Salamooni \newline
     School of Mathematical Sciences, Swami Ramanand Teerth Marathwada University, Nanded-431606, India}
      \email{ayousss83@gmail.com}

 \address{D. D. Pawar \newline
    School of Mathematical Sciences, Swami Ramanand Teerth Marathwada University, Nanded-431606, India}
     \email{dypawar@yahoo.com}


    \begin{abstract}
    In this paper, we used some theorems of fixed point for studying the
    results of
    existence and uniqueness  for Hilfer-Hadamard-Type fractional differential equations,
     \[_{H}D^{\alpha,\beta}x(t)+f(t,x(t))=0, ~~~~~~        on~~the~~ interval~~  J:=(1,e]\]
    with nonlinear boundary value problems
   \[x(1+\epsilon)=\sum_{i=1}^{n-2}\nu_{i}x(\zeta_{i}),\quad\quad\quad~_{H}D^{1,1}x(e)=\sum_{i=1}^{n-2}
\sigma_{i}~_{H}D^{1,1}x(\zeta_{i})\]
       \\\\ \textbf{AMS Classification- 34A08, 35R11}
    \end{abstract}

    \maketitle \numberwithin{equation}{section}
    \newtheorem{theorem}{Theorem}[section]
    \newtheorem{lemma}[theorem]{Lemma}
    \newtheorem{definition}[theorem]{Definition}
    \newtheorem{example}[theorem]{Example}

    \newtheorem{remark}[theorem]{Remark}
    \allowdisplaybreaks
  \[\textbf{1.Introduction}.\]
\\In this paper, we discussion the existence and uniqueness result
of the solutions for the n-point nonlinear boundary value problems
for Hilfer-Hadamard-type fractional differential equations of the
form
\begin{align*}\quad\quad\quad\quad\quad
&~_{H}D^{\alpha,\beta}x(t)+f(t,x(t))=0,\quad\quad t\in
J:=(1,e]\\&\quad\quad\quad\quad\quad\quad\quad\quad\quad\quad\quad\quad\quad\quad\quad\quad\quad\quad
\quad\quad\quad\quad\quad\quad\quad\quad\quad\quad\quad\quad(1.1)
\\&x(1+\epsilon)=\sum_{i=1}^{n-2}\nu_{i}x(\zeta_{i}),\quad\quad\quad~_{H}D^{1,1}x(e)=\sum_{i=1}^{n-2}
\sigma_{i}~_{H}D^{1,1}x(\zeta_{i})
\end{align*}
where $~_{H}D^{\alpha,\beta}$ is the Hilfer-Hadamard fractional
derivative of order $1<\alpha\leq2$ and type $\beta\in[0,1],$
$f:J\times\mathbb{R}\rightarrow\mathbb{R}$ is a continuous function,
$0<\epsilon<1,\zeta_{i}\in(1,e),\\
\nu_{i},\sigma_{i}\in\mathbb{R}\quad for\quad all\quad
i=1,2,...,n-2,\zeta_{1}<\zeta_{2}<...<\zeta_{n-2}$ and
$~_{H}D^{1,1}=t\frac{d}{dt}.$
\par The fractional differential equations give proofs of the more
 appropriate models for describing real world problems. Indeed, these
 problems cannot be described using classical integer order differential
 equations. In the past years the theory of fractional differential
 equations has received much attention from the authors, and has
 become an important field of investigation due to existence
 applications in engineering, biology, chemistry, economics and numerous branches of physics sciences[1,6,9,10].
 Fractional differential equations have a several kinds of fractional differential
 equations. One of them is the Hadamard fractional derivative innovated by
 Hadamard in 1892[4],which differs from the Riemann-Liouvill and
 Caputo type fractional derivative[9],the preceding ones in the
 sense that the kernel of the integral contains logarithmic function
 of arbitrary exponent. The properties of Hadamard fractional integral and derivative
 can be found in[1,27].Recently, the authors studied the Hadamard-type,
 Caputo-Hadamard-type and Hilfer-Hadamard-type
 fractional derivative by using the fixed point theorems with the boundary
 value problems and give the results of existence and
 uniqueness[14-21].
\par We found a variety of results for the problem (1.1) by using
traditional fixed point theorems. The first result, Theorem 3.2,
depend on Banach's Principle of contraction mapping and concerns an
existence and uniqueness result for the solutions of the problem
(1.1). In Theorem 3.3 we proved a second result of existence and
uniqueness, through a fixed point theorem and nonlinear contractions
due to Boyd and Wong. A third result of existence is proved in
Theorem 3.4, by using Krasnoselskii's fixed point theorem, and the
last result of existence, Theorem 3.5, by using Leray-Schauder type
of nonlinear alternative for single-valued maps.

\[\textbf{2.Preliminaries}\] In this section, we introduce some
notations and definitions of Hilfer-Hadamard-type fractional
calculus.
\\\textbf{\ Definition 2.1.[1,10]~} (Riemann-Liouville fractional integral).
\\ The Riemann-Liouville integral of order $~\alpha~ > 0$ of a function
$~\varphi:[1,\infty)\rightarrow\mathbb{R}~$  is defined by
$$(I^{\alpha}\varphi)(t)=\frac{1}{\Gamma(\alpha)}
_{1}\int^{t}\frac{\varphi(\tau)d\tau}{(t-\tau)^{1-\alpha}}\quad,\quad(t>1),$$
Here $\Gamma(\alpha)$ is the Euler's Gamma function.
\\\textbf{\ Definition 2.2.[1,10]~} (Riemann-Liouville fractional derivative).
\\The Riemann-Liouville fractional derivative of order $~\alpha > 0~$
of a function \\$~\varphi:[1,\infty)\rightarrow\mathbb{R}~$  is
defined by \\ \par $(D^{\alpha}\varphi)(t):=(\frac{d}{dt})^{n}
(I^{n-\alpha}\varphi)(t)$
$$\quad\quad=\frac{1}{\Gamma(n-\alpha)}\frac{d^{n}}{dt^{n}}
_{1}\int^{t}\frac{\varphi(\tau)d\tau}{(t-\tau)^{\alpha-n+1}}
\quad,\quad\quad(n=[\alpha]+1;t>1),$$ Here [$\alpha$] is the integer
part of $\alpha.$
 \\\textbf{\ Definition 2.3.[1]~}(Hadamard fractional integral).
\\~The Hadamard fractional integral of order $~\alpha\in \mathbb{R}^{+}~$for a function
$~\varphi:[1,\infty)\rightarrow\mathbb{R}~$ is defined
as\[_{H}I^{\alpha}\varphi(t)=\frac{1}{\Gamma(\alpha)}_{1}\int^{t}(\log\frac{t}{\tau})^{\alpha-1}
\quad\frac{\varphi(\tau)}{\tau}d\tau,\quad\quad(t>1)\] where
$~\log(.)=\log_{e}(.)~$.
\\\textbf{\
Definition 2.4.[1]~}(Hadamard fractional derivative).
\\~The Hadamard fractional derivative of order $~\alpha~$ applied to
the function\\ $~\varphi:[1,\infty)\rightarrow\mathbb{R}~$ is
defined
as\[_{H}D^{\alpha}\varphi(t)=\delta^{n}(_{H}I^{n-\alpha}\varphi(t)),\quad
n-1<\alpha<n,\quad n=[\alpha]+1,\]
where$\quad~\delta^{n}=(t\frac{d}{dt})^{n}\quad~$and$~[\alpha]~$denotes
the integer part of the real number$~\alpha.~$
\\\textbf{\ Definition 2.5.[4,12]~}(Caputo-Hadamard fractional derivative).
\\~The Caputo-Hadamard fractional derivative of
order$~\alpha~$applied to the function $~\varphi\in
AC_{\delta}^{n}[a,b]~$is defined
as\[_{HC}D^{\alpha}\varphi(t)=(_{H}I^{n-\alpha}\delta^{n}\varphi)(t)\]
where\\$~n=[\alpha]+1,~$and$~\varphi\in
AC_{\delta}^{n}[a,b]=\bigg\{\varphi:[a,b]\rightarrow
\mathbb{C}:\delta^{(n-1)}\varphi\in
AC[a,b],\delta=t\frac{d}{dt}\bigg\}~$
\\\textbf{\
Definition 2.6.[6,20]~}(Hilfer fractional derivative).
\\~Let$~~n-1<\alpha<n,~~0\leq\beta\leq 1,~~\varphi\in
L^{1}(a,b).~$The Hilfer fractional derivative $D^{\alpha,\beta}$ of
order$~\alpha~$and type $~\beta~$of$~\varphi~$is defined
as\[~(D^{\alpha,\beta}\varphi)(t)=\big(I^{\beta(n-\alpha)}(\frac{d}{dt})^{n}~I^{(n-\alpha)(1-\beta)}\varphi\big)(t)~\]
\[=\big(I^{\beta(n-\alpha)}(\frac{d}{dt})^{n}I^{n-\gamma}\varphi\big)(t);~\quad\gamma=\alpha+n\beta-\alpha\beta.\]
\[=\big(I^{\beta(n-\alpha)}D^{\gamma}\varphi\big)(t),\]
\\Where $I^{(.)}$ and  $~D^{(.)}~$is the Riemann-Liouvill fractional integral and
derivative defined by (2.1) and (2.2), respectively.\\In particular,
if $\quad0<\alpha<1,$ then
$$~(D^{\alpha,\beta}\varphi)(t)=\big(I^{\beta(1-\alpha)}\frac{d}{dt}~I^{(1-\alpha)(1-\beta)}\varphi\big)(t)~$$
\[\quad\quad=\big(I^{\beta(1-\alpha)}\frac{d}{dt}I^{1-\gamma}\varphi\big)(t);\quad\gamma=\alpha+\beta-\alpha\beta.~\]
\[=\big(I^{\beta(1-\alpha)}D^{\gamma}\varphi\big)(t).\]
\\ \textbf {\ Properties 2.7.[20,21].}\\ Let$~~0<\alpha<1,\quad0\leq\beta\leq 1,
\quad\gamma=\alpha+\beta-\alpha\beta,$ and $\varphi\in L^{1}(a,b).~$
If $D^{\gamma}\varphi$ exists and in $L^{1}(a,b),$ then
\[I_{a+}^{\alpha}~(D_{a+}^{\alpha,\beta}\varphi)(t)=I_{a+}^{\gamma}~
 (D_{a+}^{\gamma}\varphi)(t)=
 \varphi(t)-\frac{(I_{a+}^{1-\gamma}\varphi)(a)}{\Gamma(\gamma)}(t-a)^{\gamma-1}\]
\\\textbf{\ Definition 2.8.[20,21]}(Hilfer-Hadamard fractional derivative).
\\~Let$~~0<\alpha<1,~~0\leq\beta\leq 1,~~\varphi\in
L^{1}(a,b).~$The Hilfer-Hadamard fractional derivative
$_{H}D^{\alpha,\beta}$ of order$~\alpha~$and type
$~\beta~$of$~\varphi~$is defined
as\[~(_{H}D^{\alpha,\beta}\varphi)(t)=\big(_{H}I^{\beta(1-\alpha)}\delta~_{H}I^{(1-\alpha)(1-\beta)}\varphi\big)(t)~\]
\[=\big(_{H}I^{\beta(1-\alpha)}\delta~_{H}I^{1-\gamma}\varphi\big)(t);\quad\gamma=\alpha+\beta-\alpha\beta.\]
\[=\big(_{H}I^{\beta(1-\alpha)}_{H}D^{\gamma}\varphi\big)(t).\]
Where $_{H}I^{(.)}$ and $~_{H}D^{(.)}~$is the Hadamard fractional
integral and derivative defined by (2.3) and (2.4), respectively.
\\ \textbf {\ Theorem 2.9.[1,4].}
 \\~Let$~\Re(\alpha)>0,\quad~n=[\Re(\alpha)]+1~$and$~0<a<b<\infty.~$ if $~\varphi\in L^{1}(a,b)~$ and
 $~(_{H}I_{a+}^{n-\alpha}\varphi)(t)\in AC_{\delta}^{n}[a,b],~$ then
 \[(_{H}I_{a+}^{\alpha}~_{H}D_{a+}^{\alpha}\varphi)(t)=
 \varphi(t)-\sum_{j=0}^{n-1}\frac{(\delta^{(n-j-1)}(_{H}I_{a+}^{n-\alpha}\varphi))(a)}{\Gamma(\alpha-j)}(\log\frac{t}{a})^{\alpha-j-1}\]
\\\textbf{\ Theorem 2.10.[4,12]~}\\Let  $\varphi(t)\in AC_{\delta}^{n}[a,b]$ or $\varphi(t)\in
C_{\delta}^{n}[a,b],$ and $~~\alpha\in\mathbb{C},$
then\[(_{H}I_{a+}^{\alpha}~_{HC}D_{a+}^{\alpha}\varphi)(t)=
 \varphi(t)-\sum_{K=0}^{n-1}\frac{\delta^{K}\varphi(a)}{\Gamma(K+1)}(\log\frac{t}{a})^{K}\]
\\\textbf{\ Definition 2.11.[28] } Let $E$ be a Banach space and let $F :E\rightarrow E $ be a
mapping. $F$ is said to be a nonlinear contraction if there exists a
continuous nondecreasing function
$\Psi:\mathbb{R}^{+}\rightarrow\mathbb{R}^{+}$ such that $\Psi(0) =
0$ and $\Psi(\phi) < \phi$ for all $\phi > 0$ with the property
$$\| Fx-Fy\|\leq\Psi(\| x-y\|),\quad\quad \forall x,y\in E.$$
\\\textbf{\ Lemma 2.12.[24]} Let $E$ be a Banach space and let $F :E\rightarrow E $ be
a nonlinear contraction. Then $F$ has a unique fixed point in $E$.
\\\textbf{\ Theorem 2.13.[23]~}(Krasnoselskii's fixed point theorem). Let
$M$ be a closed, bounded, convex, and nonempty subset of a Banach
space $X.$ Let $A,B$ be the operators such that \\(a) $Ax+By\in M,$
whenever$x,y\in M;$\\(b) $A$ is compact and continuous;\\(c) $B$ is
a contraction mapping. Ten there exists $z\in M$ such that
$z=Az+Bz.$
\\\textbf{\ Theorem 2.14.[25]}(nonlinear alternative for single-valued maps). Let $E$ be a
Banach space,$ C$ a closed, convex subset of $E, U$ an open subset
of $C,$ and $0\in U.$ Suppose that $F : \overline{U}\rightarrow C
$is a continuous, compact $\big(i.e., F(\overline{U})$ is a
relatively compact subset of $C\big)$ map. Ten either
\par(i) F has a fixed point in $\overline{U}$ or
\par(ii) there is a $u\in\partial U$ (the boundary of $U$ in $C$) and $\bar{\lambda}\in (0,
1),$ with $u =\bar{\lambda} F(u).$
\\\textbf{\ Definition 2.15.[29]}(Hilfer-Hadamard fractional derivative).
\\~Let$~~n-1<\alpha<n,~~0\leq\beta\leq 1,~~\varphi\in
L^{1}(a,b).~$The Hilfer-Hadamard fractional derivative
$_{H}D^{\alpha,\beta}$ of order$~\alpha~$and type
$~\beta~$of$~\varphi~$is defined
as\[~(_{H}D^{\alpha,\beta}\varphi)(t)=\big(_{H}I^{\beta(n-\alpha)}(\delta)^{n}~_{H}I^{(n-\alpha)(1-\beta)}\varphi\big)(t)~\]
\[=\big(_{H}I^{\beta(n-\alpha)}(\delta)^{n}~_{H}I^{n-\gamma}\varphi\big)(t);\quad\gamma=\alpha+n\beta-\alpha\beta.\]
\[=\big(I^{\beta(n-\alpha)}_{H}D^{\gamma}\varphi\big)(t),\]
Where $_{H}I^{(.)}$ and $~_{H}D^{(.)}~$is the Hadamard fractional
integral and derivative defined by (2.3) and (2.4), respectively.
\\ \textbf{\ Lemma2.16.[29]}
 \\~Let$~\Re(\alpha)>0,\quad0\leq\beta\leq1,\quad\gamma=\alpha+n\beta-\alpha\beta,
 \quad n-1<\gamma\leq n,\quad~n=[\Re(\alpha)]+1~$ and $~0<a<b<\infty.~$ if $~\varphi\in L^{1}(a,b)~$ and
 $~(_{H}I_{a+}^{n-\gamma}\varphi)(t)\in AC_{\delta}^{n}[a,b],
\quad ~~$ then
 \[_{H}I_{a+}^{\alpha}~(_{H}D_{a+}^{\alpha,\beta}\varphi)(t)=_{H}I_{a+}^{\gamma}~(_{H}D_{a+}^{\gamma}\varphi)(t)=
 \varphi(t)-\sum_{j=0}^{n-1}\frac{(\delta^{(n-j-1)}(_{H}I_{a+}^{n-\gamma}\varphi))(a)}
 {\Gamma(\gamma-j)}(\log\frac{t}{a})^{\gamma-j-1}\]
 From this Lemma, we notice that if $\beta=0$ the formulae reduces to the formulae in
 the theorem 2.9, and if the $\beta=1$ the formulae reduces to the formulae in
 the theorem 2.10.\\
\\\[\textbf{3.Main Results }.\]
\\ \textbf{\ Lemma3.1}
 \\For $1<\alpha\leq2~~$,$~~0\leq\beta\leq1~~$ and $~~\varphi\in
C([1,e],\mathbb{R}),$\\$~~~\gamma=\alpha+2\beta-\alpha\beta~~$,$~\gamma\in(1,2]$
\\ the problem
\\$_{H}D^{\alpha,\beta}x(t)+\varphi(t)=0,~~~~~~~~~~~~$ $t\in J~~,$ $1<\alpha\leq2~~$,$~~0\leq\beta\leq1~~$
\\$x(1+\epsilon)=\sum_{i=1}^{n-2}\nu_{i}x(\zeta_{i}),\quad\quad\quad~$$_{H}D^{1,1}x(e)=\sum_{i=1}^{n-2}
\sigma_{i}~_{H}D^{1,1}x(\zeta_{i})\quad\quad\quad\quad\quad \quad\quad\quad\quad\quad\quad\quad\quad(3.1)$
has a unique solution it giving in the formulae
\begin{align*}
x(t)=&- _{H}I^{\alpha}\varphi(t)+\frac{(\gamma-1)\delta_{1}(\log
t)^{\gamma-2}-(\gamma-2)\delta_{2}(\log
t)^{\gamma-1}}{\lambda}\bigg[~_{H}I^{\alpha}\varphi(1+\epsilon)-\sum_{i=1}^{n-2}\nu_{i}~_{H}I^{\alpha}\varphi(\zeta_{i})\bigg]
\\&\quad+\frac{\mu_{2}(\log
t)^{\gamma-1}-\mu_{1}(\log
t)^{\gamma-2}}{\lambda}\bigg[~_{H}I^{\alpha-1}\varphi(e)-\sum_{i=1}^{n-2}\sigma_{i}~_{H}I^{\alpha-1}\varphi(\zeta_{i})\bigg],\quad\quad
t\in J.
\end{align*}
Where,
\begin{align*}
&\lambda=(\gamma-1)\delta_{1}\mu_{2}-(\gamma-2)\delta_{2}\mu_{1},\quad\quad
with \quad\quad\lambda\neq 0, \\&
\mu_{1}=(\log(1+\epsilon))^{\gamma-1}-\sum_{i=1}^{n-2}\nu_{i}(\log(\zeta_{i}))^{\gamma-1},\\&
\mu_{2}=(\log(1+\epsilon))^{\gamma-2}-\sum_{i=1}^{n-2}\nu_{i}(\log(\zeta_{i}))^{\gamma-2},\\&
\delta_{1}=1-\sum_{i=1}^{n-2}\sigma_{i}(\log(\zeta_{i}))^{\gamma-2},\\&
\delta_{2}=1-\sum_{i=1}^{n-2}\sigma_{i}(\log(\zeta_{i}))^{\gamma-3}.
\end{align*}
\\ \textbf{\ Proof.} In the view of the Lemma $(2.16)$, the solution of
the Hilfer-Hadamard differential equation $(3.1)$ can be written as
\[\quad\quad\quad\quad x(t)=-~_{H}I^{\alpha}\varphi(t)+c_{0}(\log t)^{\gamma-1}+c_{1}(\log t)^{\gamma-2}~
\quad\quad\quad\quad\quad\quad\quad\quad\quad\quad\quad~(3.2)\]
and\[\quad\quad\quad _{H}D^{1,1}
x(t)=-~_{H}I^{\alpha-1}\varphi(t)+(\gamma-1)c_{0}(\log
t)^{\gamma-2}+(\gamma-2) c_{1}(\log
t)^{\gamma-3}~\quad\quad\quad~(3.3)\] The boundary condition
$x(1+\epsilon)=\sum_{i=1}^{n-2}\nu_{i}x(\zeta_{i})$ gives
$$\quad\quad c_{1}=\frac{1}{\mu_{2}}\bigg[~_{H}I^{\alpha}\varphi(1+\epsilon)-\sum_{i=1}^{n-2}\nu_{i}~_{H}I^{\alpha}\varphi(\zeta_{i})-c_{0}\mu_{1}\bigg]
\quad\quad\quad\quad\quad\quad\quad\quad\quad\quad\quad\quad(3.4)$$
where
$$\mu_{1}=(\log(1+\epsilon))^{\gamma-1}-\sum_{i=1}^{n-2}\nu_{i}(\log(\zeta_{i}))^{\gamma-1},\quad\quad
\mu_{2}=(\log(1+\epsilon))^{\gamma-2}-\sum_{i=1}^{n-2}\nu_{i}(\log(\zeta_{i}))^{\gamma-2}.
$$
In view of the boundary condition $_{H}D^{1,1}x(e)=\sum_{i=1}^{n-2}
\sigma_{i}~_{H}D^{1,1}x(\zeta_{i})$ and by $~(3.3)~,~and~ (3.4)$ we
have
$$\quad\quad
c_{0}=\frac{1}{(\gamma-1)\delta_{1}}\bigg[-(\gamma-2)c_{1}\delta_{2}+~_{H}I^{\alpha-1}\varphi(e)-
\sum_{i=1}^{n-2}\sigma_{i}~_{H}I^{\alpha-1}\varphi(\zeta_{i})\bigg]\quad\quad\quad\quad\quad\quad(3.5)$$
where$$\delta_{1}=1-\sum_{i=1}^{n-2}\sigma_{i}(\log(\zeta_{i}))^{\gamma-2},\quad\quad
\delta_{2}=1-\sum_{i=1}^{n-2}\sigma_{i}(\log(\zeta_{i}))^{\gamma-3}.
$$
by $(3.5)$ in $(3.4)$ we have
\begin{align*}
&c_{1}=\frac{1}{\lambda}\Biggr[(\gamma-1)\delta_{1}\bigg[
~_{H}I^{\alpha}\varphi(1+\epsilon)-\sum_{i=1}^{n-2}\nu_{i}~_{H}I^{\alpha}\varphi(\zeta_{i})\bigg]\\&
\quad\quad\quad\quad\quad\quad\quad-
\mu_{1}\bigg[~_{H}I^{\alpha-1}\varphi(e)-
\sum_{i=1}^{n-2}\sigma_{i}~_{H}I^{\alpha-1}\varphi(\zeta_{i})\bigg]\Biggr]
\end{align*}
where
$$\lambda=(\gamma-1)\delta_{1}\mu_{2}-(\gamma-2)\delta_{2}\mu_{1},\quad\quad
with \quad\quad\lambda\neq 0.$$ by substituting the value of
$~c_{1}~$ in $~(3.5)~$ we have
\begin{align*}
c_{0}&=\frac{1}{\lambda}\Biggr[-(\gamma-2)\delta_{2}\bigg[
~_{H}I^{\alpha}\varphi(1+\epsilon)-\sum_{i=1}^{n-2}\nu_{i}~_{H}I^{\alpha}\varphi(\zeta_{i})\bigg]\\&
\quad\quad\quad\quad\quad\quad\quad+
\mu_{2}\bigg[~_{H}I^{\alpha-1}\varphi(e)-
\sum_{i=1}^{n-2}\sigma_{i}~_{H}I^{\alpha-1}\varphi(\zeta_{i})\bigg]\Biggr]
\end{align*}
Now substituting the values of $~c_{0}~$ and $~c_{1}~$ in $~(3.2)~$
we obtain the solution of the problem(3.1).
\[\textbf{Results of Existence}.\]
Suppose that
$$\quad\quad\quad\quad\quad\quad\quad\quad\quad\quad\quad\quad~K=C([1,e],\mathbb{R})\quad\quad\quad\quad\quad\quad
\quad\quad\quad\quad\quad\quad\quad\quad\quad\quad\quad\quad(3.6)$$
is a Banach space of all continuous functions from $[1,e]$ into
$\mathbb{R}~$ talented with the norm $~\| x\|=\sup_{t\in J}|
x(t)|.~$
\\\par From the Lemma3.1,we getting an operator $~\rho:K\rightarrow K$
defined as
\begin{align*}
\quad\quad\quad(\rho x)(t)&=-
_{H}I^{\alpha}f(\tau,x(\tau))(t)\\&\quad+\frac{(\gamma-1)\delta_{1}(\log
t)^{\gamma-2}-(\gamma-2)\delta_{2}(\log
t)^{\gamma-1}}{\lambda}\bigg[~_{H}I^{\alpha}f(\tau,x(\tau))(1+\epsilon)\\&\quad\quad\quad\quad
\quad\quad\quad\quad\quad\quad\quad-\sum_{i=1}^{n-2}\nu_{i}~_{H}I^{\alpha}f(\tau,x(\tau))(\zeta_{i})\bigg]
\quad\quad\quad\quad\quad(3.7)
\\&\quad+\frac{\mu_{2}(\log
t)^{\gamma-1}-\mu_{1}(\log
t)^{\gamma-2}}{\lambda}\bigg[~_{H}I^{\alpha-1}f(\tau,x(\tau))(e)\\&\quad\quad\quad\quad
\quad\quad\quad\quad-\sum_{i=1}^{n-2}\sigma_{i}~_{H}I^{\alpha-1}f(\tau,x(\tau))(\zeta_{i})\bigg],\quad\quad
with \quad\lambda\neq 0
\end{align*}
It must be noticed that the problem $(1.1)$ has solutions if and
only if the operator $\rho$ has fixed points.The result of existence
and uniqueness is based on the Banach Principle of contraction
mapping.
\\\textbf{\ Theorem 3.2} Let
$f:J\times\mathbb{R}\rightarrow\mathbb{R}$ be a continuous function
satisfying the supposition that \par $(Q_{1})$there exists a
constant $C > 0 $ such that $| f(t,x)-f(t,y)|\leq C| x-y|,$ for each
$t\in J$ and $x,y\in\mathbb{R}.$
\\If $~\Phi$ satisfied the condition  $~~C\Phi<1,~$ where
\begin{align*}
\quad\quad\Phi&=\bigg\{\frac{1}{\Gamma(\alpha+1)}+
\frac{(|\gamma-1|)|\delta_{1}|+(|\gamma-2|)|\delta_{2}|}{|\lambda|\Gamma(\alpha+1)}\bigg[(\log
(1+\epsilon))^{\alpha}+\sum_{i=1}^{n-2}|\nu_{i}|(\log
(\zeta_{i}))^{\alpha}\bigg]\\&\quad\quad\quad\quad+\frac{|\mu_{2}|+|\mu_{1}|}{|\lambda|\Gamma(\alpha)}
\bigg[1+\sum_{i=1}^{n-2}|\sigma_{i}|(\log(\zeta_{i}))^{\alpha-1}\bigg]
\bigg\}\quad\quad\quad\quad\quad\quad\quad\quad\quad\quad(3.8)
\end{align*}
Then the problem of boundary value $(1.1)$ has a unique solution on
$J$.
\\\textbf{\ Proof.}We using Banach's Principle of contraction
mapping for transform the problem of boundary value(1,1) into a
fixed point problem,$x=\rho x,$  where the operator $\rho$ is
defined by (3.7), we will show that $\rho$ has a fixed point which
is a unique solution of problem (1,1).\par We put $~~sup_{t\in J}|
f(\tau,0)|= p <\infty~$ and choose
\begin{align*}
\quad\quad\quad\quad\quad\quad\quad\quad\quad\quad\quad\quad\quad\quad
r\geq\frac{\Phi
P}{1-C\Phi},\quad\quad\quad\quad\quad\quad\quad\quad\quad\quad\quad\quad\quad\quad(3.9)
\end{align*}
\par Now,assume that $B_{r}=\{x\in K:| x|\leq r\},$
then we show that $\rho B_{r}\subset B_{r}.$ \\For any $x\in B_{r},$
we have
\begin{align*}
&\|\rho x\|=sup_{t\in
J}\bigg\{\Biggl|-~_{H}I^{\alpha}f(\tau,x(\tau))(t)\\&\quad\quad\quad\quad\quad\quad+\frac{(\gamma-1)\delta_{1}(\log
t)^{\gamma-2}-(\gamma-2)\delta_{2}(\log
t)^{\gamma-1}}{\lambda}\bigg[~_{H}I^{\alpha}f(\tau,x(\tau))(1+\epsilon)\\&\quad\quad\quad\quad
\quad\quad\quad\quad\quad\quad\quad\quad\quad-\sum_{i=1}^{n-2}\nu_{i}~_{H}I^{\alpha}f(\tau,x(\tau))(\zeta_{i})\bigg]
\\&\quad\quad\quad\quad\quad\quad+\frac{\mu_{2}(\log
t)^{\gamma-1}-\mu_{1}(\log
t)^{\gamma-2}}{\lambda}\bigg[~_{H}I^{\alpha-1}f(\tau,x(\tau))(e)\\&\quad\quad\quad\quad
\quad\quad\quad\quad\quad\quad\quad\quad\quad-\sum_{i=1}^{n-2}
\sigma_{i}~_{H}I^{\alpha-1}f(\tau,x(\tau))(\zeta_{i})\bigg]\Biggl|\bigg\}
\\&\quad\quad\leq sup_{t\in
J}\bigg\{~_{H}I^{\alpha}|f(\tau,x(\tau))|(t)\\&\quad\quad\quad\quad\quad\quad+\frac{(|\gamma-1|)|\delta_{1}|(\log
t)^{\gamma-2}+(|\gamma-2|)|\delta_{2}|(\log
t)^{\gamma-1}}{|\lambda|}\bigg[~_{H}I^{\alpha}|f(\tau,x(\tau))|(1+\epsilon)\\&\quad\quad\quad\quad
\quad\quad\quad\quad\quad\quad\quad\quad\quad+\sum_{i=1}^{n-2}|\nu_{i}|~_{H}I^{\alpha}|f(\tau,x(\tau))|(\zeta_{i})\bigg]
\\&\quad\quad\quad\quad\quad\quad+\frac{|\mu_{2}|(\log
t)^{\gamma-1}+|\mu_{1}|(\log
t)^{\gamma-2}}{|\lambda|}\bigg[~_{H}I^{\alpha-1}|f(\tau,x(\tau))|(e)\\&\quad\quad\quad\quad
\quad\quad\quad\quad\quad\quad\quad\quad\quad+\sum_{i=1}^{n-2}
|\sigma_{i}|~_{H}I^{\alpha-1}|f(\tau,x(\tau))|(\zeta_{i})\bigg]\bigg\}
\\&\quad\quad\leq ~_{H}I^{\alpha}\big(|f(\tau,x(\tau))-f(\tau,0)|+|f(\tau,0)|\big)(e)\\&
\quad\quad\quad\quad\quad\quad+\frac{(|\gamma-1|)|\delta_{1}|+(|\gamma-2|)|\delta_{2}|}{|\lambda|}
\bigg[~_{H}I^{\alpha}\big(|f(\tau,x(\tau))-f(\tau,0)|+|f(\tau,0)|\big)(1+\epsilon)\\&
\quad\quad\quad\quad \quad\quad\quad\quad\quad\quad\quad\quad\quad
+\sum_{i=1}^{n-2}|\nu_{i}|~_{H}I^{\alpha}\big(|f(\tau,x(\tau))-f(\tau,0)|+|f(\tau,0)|\big)(\zeta_{i})\bigg]
\\&\quad\quad\quad\quad\quad\quad+\frac{|\mu_{2}|+|\mu_{1}|}{|\lambda|}
\bigg[~_{H}I^{\alpha-1}\big(|f(\tau,x(\tau))-f(\tau,0)|+|f(\tau,0)|\big)(e)\\&\quad\quad\quad\quad
\quad\quad\quad\quad\quad\quad\quad\quad\quad+\sum_{i=1}^{n-2}
|\sigma_{i}|~_{H}I^{\alpha-1}\big(|f(\tau,x(\tau))-f(\tau,0)|+|f(\tau,0)|\big)(\zeta_{i})\bigg]
\\&\quad\quad\leq(Cr+P)\bigg\{\frac{1}{\Gamma(\alpha+1)}+
\frac{(|\gamma-1|)|\delta_{1}|+(|\gamma-2|)|\delta_{2}|}{|\lambda|\Gamma(\alpha+1)}\bigg[(\log
(1+\epsilon))^{\alpha}+\sum_{i=1}^{n-2}|\nu_{i}|(\log
(\zeta_{i}))^{\alpha}\bigg]\\&\quad\quad\quad\quad\quad\quad\quad\quad\quad+\frac{|\mu_{2}|+|\mu_{1}|}{|\lambda|\Gamma(\alpha)}
\bigg[1+\sum_{i=1}^{n-2}|\sigma_{i}|(\log(\zeta_{i}))^{\alpha-1}\bigg]
\bigg\}\\&\quad\quad=(Cr+P)\Phi\leq
r\quad\quad\quad\quad\quad\quad\quad\quad\quad\quad\quad\quad\quad\quad\quad\quad\quad\quad\quad\quad\quad\quad\quad\quad(3.10)
\end{align*}
Thus we shown $\rho B_{r}\subset B_{r}.$
\par Now,For $x,y\in K$ and
$\forall t\in J,$ we have
\begin{align*}
&\mid(\rho x)(t)-(\rho
y)(t)\mid\\&\quad=\Biggl|-~_{H}I^{\alpha}\big(f(\tau,x(\tau))-f(\tau,y(\tau))\big)(t)\\&
\quad\quad\quad\quad\quad\quad+\frac{(\gamma-1)\delta_{1}(\log
t)^{\gamma-2}-(\gamma-2)\delta_{2}(\log
t)^{\gamma-1}}{\lambda}\bigg[~_{H}I^{\alpha}\big(f(\tau,x(\tau))-f(\tau,y(\tau))\big)(1+\epsilon)\\&\quad\quad\quad\quad
\quad\quad\quad\quad\quad\quad\quad\quad\quad
-\sum_{i=1}^{n-2}\nu_{i}~_{H}I^{\alpha}\big(f(\tau,x(\tau))-f(\tau,y(\tau))\big)(\zeta_{i})\bigg]
\\&\quad\quad\quad\quad\quad\quad+\frac{\mu_{2}(\log
t)^{\gamma-1}-\mu_{1}(\log
t)^{\gamma-2}}{\lambda}\bigg[~_{H}I^{\alpha-1}\big(f(\tau,x(\tau))-f(\tau,y(\tau))\big)(e)\\&\quad\quad\quad\quad
\quad\quad\quad\quad\quad\quad\quad\quad\quad-\sum_{i=1}^{n-2}
\sigma_{i}~_{H}I^{\alpha-1}\big(f(\tau,x(\tau))-f(\tau,y(\tau))\big)(\zeta_{i})\bigg]\Biggl|
\\&\quad\quad\leq ~_{H}I^{\alpha}\big|f(\tau,x(\tau))-f(\tau,y(\tau))\big|(t)\\&
\quad\quad\quad\quad\quad\quad+\frac{(|\gamma-1|)|\delta_{1}|(\log
t)^{\gamma-2}+(|\gamma-2|)|\delta_{2}|(\log
t)^{\gamma-1}}{|\lambda|}\bigg[~_{H}I^{\alpha}\big|f(\tau,x(\tau))-f(\tau,y(\tau))\big|(1+\epsilon)\\&\quad\quad\quad\quad
\quad\quad\quad\quad\quad\quad\quad\quad\quad
+\sum_{i=1}^{n-2}|\nu_{i}|~_{H}I^{\alpha}\big|f(\tau,x(\tau))-f(\tau,y(\tau))\big|(\zeta_{i})\bigg]
\\&\quad\quad\quad\quad\quad\quad+\frac{|\mu_{2}|(\log
t)^{\gamma-1}+|\mu_{1}|(\log
t)^{\gamma-2}}{|\lambda|}\bigg[~_{H}I^{\alpha-1}\big|f(\tau,x(\tau))-f(\tau,y(\tau))\big|(e)\\&\quad\quad\quad\quad
\quad\quad\quad\quad\quad\quad\quad\quad\quad+\sum_{i=1}^{n-2}
|\sigma_{i}|~_{H}I^{\alpha-1}\big|f(\tau,x(\tau))-f(\tau,y(\tau))\big|(\zeta_{i})\bigg]
\\&\quad\quad\leq C \|x-y\|\bigg\{\frac{1}{\Gamma(\alpha+1)}+
\frac{(|\gamma-1|)|\delta_{1}|+(|\gamma-2|)|\delta_{2}|}{|\lambda|\Gamma(\alpha+1)}\bigg[(\log
(1+\epsilon))^{\alpha}+\sum_{i=1}^{n-2}|\nu_{i}|(\log
(\zeta_{i}))^{\alpha}\bigg]\\&\quad\quad\quad\quad\quad\quad\quad\quad\quad+\frac{|\mu_{2}|+|\mu_{1}|}{|\lambda|\Gamma(\alpha)}
\bigg[1+\sum_{i=1}^{n-2}|\sigma_{i}|(\log(\zeta_{i}))^{\alpha-1}\bigg]
\bigg\}\\&\quad\quad=C
\|x-y\|\Phi\quad\quad\quad\quad\quad\quad\quad\quad\quad\quad\quad\quad
\quad\quad\quad\quad\quad\quad\quad\quad\quad\quad\quad\quad\quad\quad(3.11)
\end{align*}
Therefore it shown that $\|(\rho x)(t)-(\rho y)(t)\|\leq C \Phi\|
x-y\|,$ where $C \Phi<1.$
\\Hence $\rho$ is a contraction. Thus by Banach's Principle of contraction
mapping, the problem $(1.1)$ has a uniqueness solution.
\\\textbf{\ Theorem 3.3.} Let
$f:J\times\mathbb{R}\rightarrow\mathbb{R}$ be a continuous function
satisfying the supposition that \par
$(Q_{2})~\big|f(t,x)-f(t,y)\big|\leq\varphi(t)\big(|x-y|/(P^{*}+|x-y|)\big),\quad
t\in J,\quad x,y\geq 0,$\par$\quad\quad where\quad
\varphi:J\rightarrow\mathbb{R^{+}}$ is continuous and a constant
$P^{*}$ is defined by
\begin{align*}\quad\quad\quad
P^{*}&=~_{H}I^{\alpha}\varphi(e)+\frac{(|\gamma-1|)|\delta_{1}|+(|\gamma-2|)|\delta_{2}|}{|\lambda|}
\bigg[~_{H}I^{\alpha}\varphi(1+\epsilon)+
\sum_{i=1}^{n-2}|\nu_{i}|~_{H}I^{\alpha}\varphi(\zeta_{i})\bigg]\\&\quad\quad\quad\quad\quad
+\frac{|\mu_{2}|+|\mu_{1}|}{|\lambda|}\bigg[~_{H}I^{\alpha-1}\varphi(e)+
\sum_{i=1}^{n-2}|\sigma_{i}|~_{H}I^{\alpha-1}\varphi(\zeta_{i})\bigg]\quad\quad\quad(3.12)
\end{align*}
Then the problem of boundary value $(1.1)$ has a unique solution on
$J$.
\\\textbf{\ Proof.}We have the operator $\rho:K\rightarrow K$
defined as (3.7) and we apply the definition 2.11. for that we
define a continuous nondecreasing function
$\Psi:\mathbb{R}^{+}\rightarrow\mathbb{R}^{+}$ as
\[\quad\quad\quad\quad\quad\quad\quad\quad\Psi(\phi)=\frac{P^{*}\phi}{P^{*}+\phi},\quad\quad\forall\phi\geq 0
\quad\quad\quad\quad\quad\quad\quad\quad\quad\quad\quad\quad\quad\quad(3.13)\]
Where the function $\Psi$ satisfies  $\Psi(0) = 0$ and $\Psi(\phi) <
\phi$ for all $\phi > 0.$ \par For any $x,y \in K$ and for each
$t\in J,$ we have
\begin{align*}
&\mid(\rho x)(t)-(\rho
y)(t)\mid\\&\quad=\Biggl|-~_{H}I^{\alpha}\big(f(\tau,x(\tau))-f(\tau,y(\tau))\big)(t)\\&
\quad\quad\quad\quad\quad\quad+\frac{(\gamma-1)\delta_{1}(\log
t)^{\gamma-2}-(\gamma-2)\delta_{2}(\log
t)^{\gamma-1}}{\lambda}\bigg[~_{H}I^{\alpha}\big(f(\tau,x(\tau))-f(\tau,y(\tau))\big)(1+\epsilon)\\&\quad\quad\quad\quad
\quad\quad\quad\quad\quad\quad\quad\quad\quad
-\sum_{i=1}^{n-2}\nu_{i}~_{H}I^{\alpha}\big(f(\tau,x(\tau))-f(\tau,y(\tau))\big)(\zeta_{i})\bigg]
\\&\quad\quad\quad\quad\quad\quad+\frac{\mu_{2}(\log
t)^{\gamma-1}-\mu_{1}(\log
t)^{\gamma-2}}{\lambda}\bigg[~_{H}I^{\alpha-1}\big(f(\tau,x(\tau))-f(\tau,y(\tau))\big)(e)\\&\quad\quad\quad\quad
\quad\quad\quad\quad\quad\quad\quad\quad\quad-\sum_{i=1}^{n-2}
\sigma_{i}~_{H}I^{\alpha-1}\big(f(\tau,x(\tau))-f(\tau,y(\tau))\big)(\zeta_{i})\bigg]\Biggl|
\\&\quad\quad\leq ~_{H}I^{\alpha}\big|f(\tau,x(\tau))-f(\tau,y(\tau))\big|(t)\\&
\quad\quad\quad\quad\quad\quad+\frac{(|\gamma-1|)|\delta_{1}|(\log
t)^{\gamma-2}+(|\gamma-2|)|\delta_{2}|(\log
t)^{\gamma-1}}{|\lambda|}\bigg[~_{H}I^{\alpha}\big|f(\tau,x(\tau))-f(\tau,y(\tau))\big|(1+\epsilon)\\&\quad\quad\quad\quad
\quad\quad\quad\quad\quad\quad\quad\quad\quad
+\sum_{i=1}^{n-2}|\nu_{i}|~_{H}I^{\alpha}\big|f(\tau,x(\tau))-f(\tau,y(\tau))\big|(\zeta_{i})\bigg]
\\&\quad\quad\quad\quad\quad\quad+\frac{|\mu_{2}|(\log
t)^{\gamma-1}+|\mu_{1}|(\log
t)^{\gamma-2}}{|\lambda|}\bigg[~_{H}I^{\alpha-1}\big|f(\tau,x(\tau))-f(\tau,y(\tau))\big|(e)\\&\quad\quad\quad\quad
\quad\quad\quad\quad\quad\quad\quad\quad\quad+\sum_{i=1}^{n-2}
|\sigma_{i}|~_{H}I^{\alpha-1}\big|f(\tau,x(\tau))-f(\tau,y(\tau))\big|(\zeta_{i})\bigg]
\\&\quad\quad\leq ~_{H}I^{\alpha}\bigg(\varphi(\tau)\frac{|x(\tau)-y(\tau)|}{P^{*}+|x(\tau)-y(\tau)|}\bigg)(e)\\&
\quad\quad\quad\quad\quad\quad
+\frac{(|\gamma-1|)|\delta_{1}|+(|\gamma-2|)|\delta_{2}|}{|\lambda|}\bigg[
~_{H}I^{\alpha}\bigg(\varphi(\tau)\frac{|x(\tau)-y(\tau)|}{P^{*}+|x(\tau)-y(\tau)|}\bigg)(1+\epsilon)\\&\quad\quad\quad\quad
\quad\quad\quad\quad\quad\quad\quad\quad\quad+\sum_{i=1}^{n-2}|\nu_{i}|~_{H}I^{\alpha}
\bigg(\varphi(\tau)\frac{|x(\tau)-y(\tau)|}{P^{*}+|x(\tau)-y(\tau)|}\bigg)(\zeta_{i})\bigg]
\\&\quad\quad\quad\quad\quad\quad+\frac{|\mu_{2}|+|\mu_{1}|}{|\lambda|}
\bigg[~_{H}I^{\alpha-1}\bigg(\varphi(\tau)\frac{|x(\tau)-y(\tau)|}{P^{*}+|x(\tau)-y(\tau)|}\bigg)(e)\\&\quad\quad\quad\quad
\quad\quad\quad\quad\quad\quad\quad\quad\quad+\sum_{i=1}^{n-2}
|\sigma_{i}|~_{H}I^{\alpha-1}\bigg(\varphi(\tau)\frac{|x(\tau)-y(\tau)|}{P^{*}+|x(\tau)-y(\tau)|}\bigg)(\zeta_{i})\bigg]
\\&\quad\quad\leq\frac{\Psi(\|x-y\|)}{P^{*}}\biggr\{~_{H}I^{\alpha}\varphi(e)+
\frac{(|\gamma-1|)|\delta_{1}|+(|\gamma-2|)|\delta_{2}|}{|\lambda|}\bigg[~_{H}I^{\alpha}\varphi(1+\epsilon)+
\sum_{i=1}^{n-2}|\nu_{i}|~_{H}I^{\alpha}\varphi(\zeta_{i})\bigg]\\&\quad\quad\quad\quad\quad
+\frac{|\mu_{2}|+|\mu_{1}|}{|\lambda|}\bigg[~_{H}I^{\alpha-1}\varphi(e)+
\sum_{i=1}^{n-2}|\sigma_{i}|~_{H}I^{\alpha-1}\varphi(\zeta_{i})\bigg]\biggr\}
\\&\quad\quad=\Psi(\|x-y\|).\quad\quad\quad\quad\quad\quad\quad\quad\quad\quad\quad\quad\quad\quad
\quad\quad\quad\quad\quad\quad\quad\quad\quad\quad\quad\quad(3.14)
\end{align*}
Which implies that $\|\rho x-\rho y\|\leq\Psi(\|x-y\|).$ Then $\rho$
is a nonlinear contraction.Thus, by Lemma 2.12 the operator $\rho$
has a fixed point, that is the unique solution of the problem
$(1.1)$ .
\par Next, we give a result of existence by using Theorem 2.13(Krasnoselskii's fixed point theorem).
\\\textbf{\ Theorem 3.4.} Let
$f:J\times\mathbb{R}\rightarrow\mathbb{R}$ be a continuous function
satisfying the supposition that $(Q_{1}).$ In addition assume that
$$(Q_{3})\quad |f(t,x)|\leq g(t),\quad \forall (t,x)\in
J\times\mathbb{R}\quad and \quad g\in C([1,e],\mathbb{R}^{+}).$$ If
$$\quad\quad\quad\quad\quad\quad\quad\quad\quad\quad\frac{C}{\Gamma(\alpha+1)}<1,\quad\quad\quad\quad\quad
\quad\quad\quad\quad\quad\quad\quad\quad\quad\quad\quad\quad\quad(3.15)$$
Then the problem of boundary value $(1.1)$ has at least one solution
on $J$.
\\\textbf{\ Proof.} We put $\sup_{t\in J}|g(t)|=\|g\|$ and choose a
suitable constant $\hat{r}$ as
\[\quad\quad\quad\quad\quad\quad\quad\quad\quad\quad\quad\quad\hat{r}\geq\|g\|\Phi,\quad\quad\quad\quad\quad
\quad\quad\quad\quad\quad\quad\quad\quad\quad\quad\quad\quad\quad(3.16)\]
where $\Phi$ is defined by (3.8).Moreover, we set the operators
$\mathscr{F}$ and $\mathscr{G}$ on \\ $B_{\hat{r}}=\{x\in K
:\|x\|\leq\hat{r}\}$ as
\begin{align*}
&(\mathscr{F}x)(t)=\frac{(\gamma-1)\delta_{1}(\log
t)^{\gamma-2}-(\gamma-2)\delta_{2}(\log
t)^{\gamma-1}}{\lambda}\bigg[~_{H}I^{\alpha}f(\tau,x(\tau))(1+\epsilon)\\&\quad\quad\quad\quad
\quad\quad\quad\quad\quad\quad\quad\quad\quad\quad-\sum_{i=1}^{n-2}\nu_{i}~_{H}I^{\alpha}f(\tau,x(\tau))(\zeta_{i})\bigg]
\quad\quad\quad\quad\quad\quad\quad\quad(3.17)
\\&\quad+\frac{\mu_{2}(\log
t)^{\gamma-1}-\mu_{1}(\log
t)^{\gamma-2}}{\lambda}\bigg[~_{H}I^{\alpha-1}f(\tau,x(\tau))(e)\\&\quad\quad\quad\quad
\quad\quad\quad\quad-\sum_{i=1}^{n-2}\sigma_{i}~_{H}I^{\alpha-1}f(\tau,x(\tau))(\zeta_{i})\bigg],\quad\quad
\quad t\in J
\\&(\mathscr{G}x)(t)= -
_{H}I^{\alpha}f(\tau,x(\tau))(t),\quad t\in J
\end{align*}
For any $x,y\in B_{\hat{r}},$ we have
\begin{align*}
\|\mathscr{F}x+\mathscr{G}x\|&\leq
\|g\|\bigg(\frac{1}{\Gamma(\alpha+1)}+
\frac{(|\gamma-1|)|\delta_{1}|+(|\gamma-2|)|\delta_{2}|}{|\lambda|\Gamma(\alpha+1)}\bigg[(\log
(1+\epsilon))^{\alpha}+\sum_{i=1}^{n-2}|\nu_{i}|(\log
(\zeta_{i}))^{\alpha}\bigg]\\&\quad\quad\quad\quad+\frac{|\mu_{2}|+|\mu_{1}|}{|\lambda|\Gamma(\alpha)}
\bigg[1+\sum_{i=1}^{n-2}|\sigma_{i}|(\log(\zeta_{i}))^{\alpha-1}\bigg]\bigg)\\&=\|g\|\Phi\leq\hat{r}.
\quad\quad\quad\quad\quad\quad\quad\quad\quad\quad\quad\quad\quad\quad
\quad\quad\quad\quad\quad\quad\quad\quad\quad(3.18)
\end{align*}
This implies $\mathscr{F}x+\mathscr{G}x\in B_{\hat{r}}.$ It follows
from supposition $(Q_{1})$ together with (3.15) that $\mathscr{G}$
is a mapping of contraction.Furthermore,it is easy to show that the
operator $\mathscr{F}$ is continuous.
\begin{align*}
\quad\quad \quad
\|\mathscr{F}x\|&\leq\|g\|\bigg(\frac{(|\gamma-1|)|\delta_{1}|+(|\gamma-2|)|\delta_{2}|}{|\lambda|\Gamma(\alpha+1)}\bigg[(\log
(1+\epsilon))^{\alpha}+\sum_{i=1}^{n-2}|\nu_{i}|(\log
(\zeta_{i}))^{\alpha}\bigg]\\&\quad\quad\quad\quad\quad+\frac{|\mu_{2}|+|\mu_{1}|}{|\lambda|\Gamma(\alpha)}
\bigg[1+\sum_{i=1}^{n-2}|\sigma_{i}|(\log(\zeta_{i}))^{\alpha-1}\bigg]\bigg).
\quad\quad \quad\quad\quad\quad(3.19)
\end{align*}
Hence,$\mathscr{F}$ is uniformly bounded on $B_{\hat{r}}.$ \par
Next,we prove that the operator $\mathscr{F}$ is a compactness, for
that we put \\$\sup_{(t,x)\in J\times
B_{\hat{r}}}|f(t,x)|=\bar{p}<\infty .$ Consequently, for
$t_{1},t_{2}\in J,$ we get
\begin{align*}
&|(\mathscr{F}x)(t_{1})-(\mathscr{F}x)(t_{2})|\\
&\quad=\Biggl|\Bigg\{\frac{(\gamma-1)\delta_{1}(\log
t_{1})^{\gamma-2}-(\gamma-2)\delta_{2}(\log
t_{1})^{\gamma-1}}{\lambda}\bigg[~_{H}I^{\alpha}f(\tau,x(\tau))(1+\epsilon)\\&\quad\quad\quad\quad\quad\quad\quad\quad
\quad\quad\quad\quad\quad\quad-\sum_{i=1}^{n-2}\nu_{i}~_{H}I^{\alpha}f(\tau,x(\tau))(\zeta_{i})\bigg]
\\&\quad\quad\quad\quad+\frac{\mu_{2}(\log
t_{1})^{\gamma-1}-\mu_{1}(\log
t_{1})^{\gamma-2}}{\lambda}\bigg[~_{H}I^{\alpha-1}f(\tau,x(\tau))(e)\\&\quad\quad\quad\quad
\quad\quad\quad\quad-\sum_{i=1}^{n-2}\sigma_{i}~_{H}I^{\alpha-1}f(\tau,x(\tau))(\zeta_{i})\bigg]\Bigg\}\\&
\quad\quad\quad-\Bigg\{\frac{(\gamma-1)\delta_{1}(\log
t_{2})^{\gamma-2}-(\gamma-2)\delta_{2}(\log
t_{2})^{\gamma-1}}{\lambda}\bigg[~_{H}I^{\alpha}f(\tau,x(\tau))(1+\epsilon)\\&\quad\quad\quad\quad\quad\quad\quad\quad
\quad\quad\quad\quad\quad\quad\quad-\sum_{i=1}^{n-2}\nu_{i}~_{H}I^{\alpha}f(\tau,x(\tau))(\zeta_{i})\bigg]
\\&\quad\quad\quad\quad+\frac{\mu_{2}(\log
t_{2})^{\gamma-1}-\mu_{1}(\log
t_{2})^{\gamma-2}}{\lambda}\bigg[~_{H}I^{\alpha-1}f(\tau,x(\tau))(e)\\&\quad\quad\quad\quad
\quad\quad\quad\quad-\sum_{i=1}^{n-2}\sigma_{i}~_{H}I^{\alpha-1}f(\tau,x(\tau))(\zeta_{i})\bigg]\Bigg\}\Biggl|\\&\quad
\leq\bar{p}\quad\frac{(|\gamma-1|)|\delta_{1}||(\log
t_{2})^{\gamma-2}-\log
t_{1})^{\gamma-2}|+(|\gamma-2|)|\delta_{2}||(\log
t_{2})^{\gamma-1}-(\log
t_{1})^{\gamma-1}|}{|\lambda|\Gamma(\alpha+1)}\bigg[(\log
(1+\epsilon))^{\alpha}\\&\quad\quad\quad\quad\quad\quad\quad\quad\quad\quad\quad+\sum_{i=1}^{n-2}|\nu_{i}|(\log
(\zeta_{i}))^{\alpha}\bigg]\\&\quad\quad+\bar{p}\quad\frac{|\mu_{2}||(\log
t_{2})^{\gamma-1}-\log t_{1})^{\gamma-1}|+|\mu_{1}||(\log
t_{2})^{\gamma-2}-(\log
t_{1})^{\gamma-2}|}{|\lambda|\Gamma(\alpha)}\bigg[1+\sum_{i=1}^{n-2}|\sigma_{i}|(\log(\zeta_{i}))^{\alpha-1}\bigg]
\\&\quad\quad\quad\quad\quad\quad\quad\quad\quad\quad\quad\quad\quad
\quad\quad\quad\quad\quad\quad\quad\quad\quad\quad\quad\quad
\quad\quad\quad\quad\quad\quad\quad\quad\quad(3.20)
\end{align*}
which is independent of $x$ and tends to zero as $t_{2}\rightarrow
t_{1}.$ Thus,$\mathscr{F}$ is equicontinuous. Hence $\mathscr{F}$ is
relatively compact on $ B_{\hat{r}}.$ Therefore, by the
Arzela-Ascoli theorem, $\mathscr{F}$ is compact on $ B_{\hat{r}}.$
Thus, by Theorem 2.13 the problem of boundary value(1.1) has at
least one solution on $J.$
\par Now, the finally result of existence
is based on Theorem 2.14(nonlinear alternative for single-valued
maps).
\\\textbf{\ Theorem 3.5.} Let
$f:J\times\mathbb{R}\rightarrow\mathbb{R}$ be a continuous function,
and assume that \par $(Q_{4})$ there exists a continuous
nondecreasing function $ \vartheta: \mathbb{R^{+}}\rightarrow
\mathbb{R^{+}}\backslash\{0\}$ such that
\[\quad\quad\quad\quad|f(t,x)|\leq q(t)\vartheta(|x|)\quad for \quad each \quad (t,x)\in J\times\mathbb{R}
\quad\quad\quad\quad\quad\quad\quad\quad(3.21)\] where $q\in
C([1,e],\mathbb{R}^{+})$ be a function.
\par $(Q_{5})$ there exists a constant $L>0$ such that
$$\quad\quad\quad\quad\quad\quad\quad\quad\quad\quad\quad\quad\quad\quad\frac{L}{\|q\|\vartheta(L)\Phi}>1,
\quad\quad\quad\quad\quad\quad\quad\quad\quad\quad\quad\quad\quad(3.22)$$
where $\Phi$ is defined by (3.8). \\Then the problem of boundary
value $(1.1)$ has at least one solution on $J$.
\\\textbf{\ Proof.} We have the operator $\rho$ is
defined by (3.7). Firstly, we will show that $\rho$ maps bounded
sets (balls) into bounded sets in $K$, for that let $\bar{r}$ a
positive number, and $B_{\bar{r}}=\{x\in K :\|x\|\leq\bar{r}\}$ be a
bounded ball in $K,$ where $K$ is defined  by (3.6). For $t\in J,$
we have
\begin{align*}
|\rho
x(t)|&\leq~_{H}I^{\alpha}|f(\tau,x(\tau))|(e)\\&\quad\quad\quad+
\frac{(|\gamma-1|)|\delta_{1}|+(|\gamma-2|)|\delta_{2}|}{|\lambda|}\bigg[~_{H}I^{\alpha}|f(\tau,x(\tau))|(1+\epsilon)\\&
\quad\quad\quad\quad\quad\quad\quad\quad+
\sum_{i=1}^{n-2}|\nu_{i}|~_{H}I^{\alpha}|f(\tau,x(\tau))|(\zeta_{i})\bigg]
\\&\quad\quad\quad+\frac{|\mu_{2}|+|\mu_{1}|}{|\lambda|}\bigg[~_{H}I^{\alpha-1}|f(\tau,x(\tau))|(e)\\&\quad\quad\quad\quad
\quad\quad\quad\quad+\sum_{i=1}^{n-2}
|\sigma_{i}|~_{H}I^{\alpha-1}|f(\tau,x(\tau))|(\zeta_{i})\bigg]\\&\leq
\|q\|\vartheta(\|x\|)\frac{1}{\Gamma(\alpha+1)}\\&\quad\quad+\|q\|\vartheta(\|x\|)
\frac{(|\gamma-1|)|\delta_{1}|+(|\gamma-2|)|\delta_{2}|}{|\lambda|\Gamma(\alpha+1)}
\bigg[(\log(1+\epsilon))^{\alpha}+
\sum_{i=1}^{n-2}|\nu_{i}|(\log(\zeta_{i}))^{\alpha}\bigg]
\\&\quad\quad\quad\quad+\|q\|\vartheta(\|x\|)\frac{|\mu_{2}|+|\mu_{1}|}{|\lambda|\Gamma(\alpha)}
\bigg[1+\sum_{i=1}^{n-2}
|\sigma_{i}|(\log(\zeta_{i}))^{\alpha-1}\bigg]\\&\leq
\|q\|\vartheta(\bar{r})\bigg\{\frac{1}{\Gamma(\alpha+1)}+
\frac{(|\gamma-1|)|\delta_{1}|+(|\gamma-2|)|\delta_{2}|}{|\lambda|\Gamma(\alpha+1)}\bigg[(\log
(1+\epsilon))^{\alpha}+\sum_{i=1}^{n-2}|\nu_{i}|(\log
(\zeta_{i}))^{\alpha}\bigg]\\&\quad\quad\quad\quad\quad\quad+\frac{|\mu_{2}|+|\mu_{1}|}{|\lambda|\Gamma(\alpha)}
\bigg[1+\sum_{i=1}^{n-2}|\sigma_{i}|(\log(\zeta_{i}))^{\alpha-1}\bigg]
\bigg\}\\&:=C_{1}.\quad\quad\quad\quad\quad\quad\quad\quad\quad\quad\quad\quad\quad
\quad\quad\quad\quad\quad\quad\quad\quad\quad\quad\quad\quad\quad\quad\quad\quad(3.23)
\end{align*}
This implies that $\|\rho x\|\leq C_{1}.$
\par Now, we will show that $\rho$ maps bounded
sets into equicontinuous sets of $K$, for that let $\sup_{(t,x)\in
J\times B_{\bar{r}}}|f(t,x)|=p^{\star}<\infty ,$
$\omega_{1},\omega_{2}\in J,$ with $\omega_{1}<\omega_{2}$  and
$x\in B_{\bar{r}}.$ Hence we have
\begin{align*}
&|(\rho x)(\omega_{1})-(\rho x)(\omega_{2})|\\
&\quad=\Biggl|\Bigg\{-
_{H}I^{\alpha}f(\tau,x(\tau))(\omega_{1})+\frac{(\gamma-1)\delta_{1}(\log
\omega_{1})^{\gamma-2}-(\gamma-2)\delta_{2}(\log
\omega_{1})^{\gamma-1}}{\lambda}\bigg[~_{H}I^{\alpha}f(\tau,x(\tau))(1+\epsilon)\\&\quad\quad\quad\quad\quad\quad\quad\quad
\quad\quad\quad\quad\quad\quad-\sum_{i=1}^{n-2}\nu_{i}~_{H}I^{\alpha}f(\tau,x(\tau))(\zeta_{i})\bigg]
\\&\quad\quad\quad\quad+\frac{\mu_{2}(\log
\omega_{1})^{\gamma-1}-\mu_{1}(\log
\omega_{1})^{\gamma-2}}{\lambda}\bigg[~_{H}I^{\alpha-1}f(\tau,x(\tau))(e)\\&\quad\quad\quad\quad
\quad\quad\quad\quad-\sum_{i=1}^{n-2}\sigma_{i}~_{H}I^{\alpha-1}f(\tau,x(\tau))(\zeta_{i})\bigg]\Bigg\}\\&
\quad\quad\quad-\Bigg\{-
_{H}I^{\alpha}f(\tau,x(\tau))(\omega_{2})+\frac{(\gamma-1)\delta_{1}(\log
\omega_{2})^{\gamma-2}-(\gamma-2)\delta_{2}(\log
\omega_{2})^{\gamma-1}}{\lambda}\bigg[~_{H}I^{\alpha}f(\tau,x(\tau))(1+\epsilon)\\&\quad\quad\quad\quad\quad\quad\quad\quad
\quad\quad\quad\quad\quad\quad\quad-\sum_{i=1}^{n-2}\nu_{i}~_{H}I^{\alpha}f(\tau,x(\tau))(\zeta_{i})\bigg]
\\&\quad\quad\quad\quad+\frac{\mu_{2}(\log
t_{2})^{\gamma-1}-\mu_{1}(\log
t_{2})^{\gamma-2}}{\lambda}\bigg[~_{H}I^{\alpha-1}f(\tau,x(\tau))(e)\\&\quad\quad\quad\quad
\quad\quad\quad\quad-\sum_{i=1}^{n-2}\sigma_{i}~_{H}I^{\alpha-1}f(\tau,x(\tau))(\zeta_{i})\bigg]\Bigg\}\Biggl|\\&\quad
\leq p^{\star}\quad\frac{|(\log \omega_{2})^{\alpha}-\log
\omega_{1})^{\alpha}|}{\Gamma(\alpha+1)}\\&\quad\quad\quad+
p^{\star}\quad\frac{(|\gamma-1|)|\delta_{1}||(\log
\omega_{2})^{\gamma-2}-\log
\omega_{1})^{\gamma-2}|+(|\gamma-2|)|\delta_{2}||(\log
\omega_{2})^{\gamma-1}-(\log
\omega_{1})^{\gamma-1}|}{|\lambda|\Gamma(\alpha+1)}\bigg[(\log
(1+\epsilon))^{\alpha}\\&\quad\quad\quad\quad\quad\quad\quad\quad\quad\quad\quad+\sum_{i=1}^{n-2}|\nu_{i}|(\log
(\zeta_{i}))^{\alpha}\bigg]\\&\quad\quad+p^{\star}\quad\frac{|\mu_{2}||(\log
\omega_{2})^{\gamma-1}-\log \omega_{1})^{\gamma-1}|+|\mu_{1}||(\log
\omega_{2})^{\gamma-2}-(\log
\omega_{1})^{\gamma-2}|}{|\lambda|\Gamma(\alpha)}\bigg[1+\sum_{i=1}^{n-2}|\sigma_{i}|(\log(\zeta_{i}))^{\alpha-1}\bigg]
\\&\quad\quad\quad\quad\quad\quad\quad\quad\quad\quad\quad\quad\quad
\quad\quad\quad\quad\quad\quad\quad\quad\quad\quad\quad\quad
\quad\quad\quad\quad\quad\quad\quad\quad\quad(3.24)
\end{align*}
Clearly, as $\omega_{2}\rightarrow\omega_{1}$ the right hand side of
the previous inequality tends to zero which is independently of
$x\in B_{\bar{r}}.$ Thus, by the Arzela-Ascoli theorem, it follows
that $\rho : K\rightarrow K$ is completely continuous.
\par Finally, let $x$ be a solution. So, for $t\in J,$ following the
similar computations as in the first step, we have
\begin{align*}
\|x\|&\leq\|q\|\vartheta(\|x\|)\frac{1}{\Gamma(\alpha+1)}\\&\quad\quad+\|q\|\vartheta(\|x\|)
\frac{(|\gamma-1|)|\delta_{1}|+(|\gamma-2|)|\delta_{2}|}{|\lambda|\Gamma(\alpha+1)}
\bigg[(\log(1+\epsilon))^{\alpha}+
\sum_{i=1}^{n-2}|\nu_{i}|(\log(\zeta_{i}))^{\alpha}\bigg]
\\&\quad\quad\quad\quad+\|q\|\vartheta(\|x\|)\frac{|\mu_{2}|+|\mu_{1}|}{|\lambda|\Gamma(\alpha)}
\bigg[1+\sum_{i=1}^{n-2}
|\sigma_{i}|(\log(\zeta_{i}))^{\alpha-1}\bigg]\\&=\|q\|\vartheta(\|x\|)\Phi.\quad\quad\quad\quad\quad\quad\quad\quad\quad\quad
\quad\quad\quad\quad\quad\quad\quad\quad\quad\quad\quad\quad\quad\quad\quad\quad(3.25)
\end{align*}
Thus, we have
$$\quad\quad\quad\quad\quad\quad\quad\quad\quad\quad\quad\quad\quad\quad\frac{\|x\|}{\|q\|\vartheta(\|x\|)\Phi}\leq 1.
\quad\quad\quad\quad\quad\quad\quad\quad\quad\quad\quad\quad\quad(3.26)$$
In view of $(Q_{5}),$ there exists $L$ such that $\| x\|\neq L.$ Let
us set
$$\quad\quad\quad\quad\quad\quad\quad\quad\quad\quad\quad
V=\{x \in K:\|x\|<L\}.
\quad\quad\quad\quad\quad\quad\quad\quad\quad\quad\quad\quad\quad(3.27)$$
Note that the operator $\rho:\overline{V}\rightarrow K$ is
continuous and completely continuous. From the choice of $V,$ there
is no $x\in\partial V$ such that $x=\bar{\lambda}\rho x$ for some
$\bar{\lambda}\in(0, 1).$ Thus, by Theorem 2.14 the operator $\rho$
has a fixed point in $\overline{V}$ which is a solution of the
problem of boundary value(1.1).
\[ \textbf{\ 4.Examples}\]
\\ \textbf{\ Example 4.1.} \\ Consider the following boundary value
problem for Hilfer-Hadamard-type fractional differential equation:
$$\quad\quad\quad\quad\quad~_{H}D^{3/2,1/2}x(t)+f(t,x(t))=0,~\quad\quad~~t\in
J=(1,e]~\quad\quad\quad\quad\quad\quad\quad\quad\quad\quad(4.1)$$
$$\quad x(1.3)=\frac{1}{2}x(3/2)-\frac{3}{4}x(7/4),~\quad\quad
~_{H}D^{1,1}x(e)=\frac{2}{3}~_{H}D^{1,1}x(3/2)+\frac{4}{3}~_{H}D^{1,1}x(7/4).~$$
Here,
$$\alpha=3/2,\quad\beta=1/2,\quad\gamma=7/4,\quad\nu_{1}=1/2,\quad\nu_{2}=-3/4
,\quad\sigma_{1}=2/3,\quad\sigma_{2}=4/3 ,
$$ $\zeta_{1}=3/2,\quad\zeta_{2}=7/4,\quad\epsilon=0.3,\quad1+\epsilon=1.3
$ and $f(t,x(t))=\frac{(\sqrt{t}+\log
t^{2})}{2e^{t}(3+t)^{2}}\big(\frac{| x(t)|}{2+| x(t)|}\big).$
\\Clearly,
$$\quad\quad\quad\quad\quad\quad\quad\quad\quad|f(t,x)-f(t,y)|\leq\frac{3}{64e}(|x-y|)
\quad\quad\quad\quad\quad\quad\quad\quad\quad\quad\quad\quad(4.2)$$
Hence$(Q_{1})$ is satisfied with $C=\frac{3}{64e}.$ We can show that
\begin{align*}
&\quad\quad\mu_{1}=(\log(1+\epsilon))^{\gamma-1}-\sum_{i=1}^{n-2}\nu_{i}(\log(\zeta_{i}))^{\gamma-1}\approx0.59779,\\&\quad\quad
\mu_{2}=(\log(1+\epsilon))^{\gamma-2}-\sum_{i=1}^{n-2}\nu_{i}(\log(\zeta_{i}))^{\gamma-2}\approx1.63780,\\&\quad\quad
\delta_{1}=1-\sum_{i=1}^{n-2}\sigma_{i}(\log(\zeta_{i}))^{\gamma-2}\approx-1.37703,\\&\quad\quad
\delta_{2}=1-\sum_{i=1}^{n-2}\sigma_{i}(\log(\zeta_{i}))^{\gamma-3}\approx-3.81518,
\\&\quad\quad\lambda=(\gamma-1)\delta_{1}\mu_{2}-(\gamma-2)\delta_{2}\mu_{1}\approx-2.26164,\\
&\quad\quad\Phi=\frac{1}{\Gamma(\alpha+1)}+
\frac{(|\gamma-1|)|\delta_{1}|+(|\gamma-2|)|\delta_{2}|}{|\lambda|\Gamma(\alpha+1)}\bigg[(\log
(1+\epsilon))^{\alpha}+\sum_{i=1}^{n-2}|\nu_{i}|(\log
(\zeta_{i}))^{\alpha}\bigg]\\&\quad\quad\quad\quad+\frac{|\mu_{2}|+|\mu_{1}|}{|\lambda|\Gamma(\alpha)}
\bigg[1+\sum_{i=1}^{n-2}|\sigma_{i}|(\log(\zeta_{i}))^{\alpha-1}\bigg]\\&\quad\quad\approx3.835201,\\&
C\Phi=\frac{3}{64e}(3.835201)\approx0.06613554378<1.
\end{align*}Therefore,
by Theorem3.2, the boundary value problem (4.1) has a unique
solution on $J.$
\\ \textbf{\ Example 4.2.} \\ Consider the following boundary value
problem for Hilfer-Hadamard-type fractional differential equation:
$$\quad\quad\quad\quad\quad~_{H}D^{3/2,2/3}x(t)+f(t,x(t))=0,~\quad\quad~~t\in
J=(1,e]~\quad\quad\quad\quad\quad\quad\quad\quad\quad\quad$$
$$\quad\quad\quad\quad\quad\quad
x(1.5)=2x(4/3)-\frac{1}{2}x(2)+\frac{5}{3}x(9/7),~\quad\quad\quad\quad\quad\quad\quad\quad\quad\quad\quad\quad(4.3)$$
$$\quad\quad\quad\quad\quad\quad~_{H}D^{1,1}x(e)=-~_{H}D^{1,1}x(4/3)+3D^{1,1}x(2)-\frac{11}{3}~_{H}D^{1,1}x(9/7).
\quad\quad\quad\quad\quad\quad\quad\quad\quad\quad\quad\quad~$$
Here,
$$\alpha=3/2,\quad\beta=2/3,\quad\gamma=11/6,\quad\nu_{1}=2,\quad\nu_{2}=-1/2,\quad\nu_{3}=5/3
,\quad\sigma_{1}=-1,\quad\sigma_{2}=3,
$$
$\sigma_{3}=-11/3
,\zeta_{1}=4/3,\quad\zeta_{2}=2,\quad\zeta_{2}=9/7,\quad\epsilon=0.5,\quad1+\epsilon=1.5
$ \\and
$$\quad\quad\quad\quad\quad\quad\quad\quad\quad f(t,x(t))=\frac{(1+\log
t)}{(t+1)^{2}}\big(\frac{| x(t)|+1}{3+|
x(t)|}\big).\quad\quad\quad\quad\quad\quad\quad\quad\quad\quad(4.4)$$
\\Clearly,
$$\quad\quad\quad\quad\quad\quad\quad\quad\quad|f(t,x)|\leq|\frac{(1+\log t)}{(t+1)^{2}}\big(\frac{|
x(t)|+1}{3+| x(t)|}\big)|
\quad\quad\quad\quad\quad\quad\quad\quad\quad\quad\quad\quad$$
$$\quad\quad\quad\quad\quad\quad\quad\quad\quad\quad\quad\quad\quad\leq|(1+\log t)\big(\frac{|
x(t)|+1}{12}\big)
\quad\quad\quad\quad\quad\quad\quad\quad\quad\quad\quad(4.5)$$ we
choose $q(t)=1+\log t$ and $\vartheta(|x|)=(| x(t)|+1)/12,$ We can
show that
\begin{align*}
&\quad\quad\mu_{1}=(\log(1+\epsilon))^{\gamma-1}-\sum_{i=1}^{n-2}\nu_{i}(\log(\zeta_{i}))^{\gamma-1}\approx-0.395713,\\&\quad\quad
\mu_{2}=(\log(1+\epsilon))^{\gamma-2}-\sum_{i=1}^{n-2}\nu_{i}(\log(\zeta_{i}))^{\gamma-2}\approx-2.865742,\\&\quad\quad
\delta_{1}=1-\sum_{i=1}^{n-2}\sigma_{i}(\log(\zeta_{i}))^{\gamma-2}\approx3.65750,\\&\quad\quad
\delta_{2}=1-\sum_{i=1}^{n-2}\sigma_{i}(\log(\zeta_{i}))^{\gamma-3}\approx19.04369,
\\&\quad\quad\lambda=(\gamma-1)\delta_{1}\mu_{2}-(\gamma-2)\delta_{2}\mu_{1}\approx-9.990516,\\
&\quad\quad\Phi=\frac{1}{\Gamma(\alpha+1)}+
\frac{(|\gamma-1|)|\delta_{1}|+(|\gamma-2|)|\delta_{2}|}{|\lambda|\Gamma(\alpha+1)}\bigg[(\log
(1+\epsilon))^{\alpha}+\sum_{i=1}^{n-2}|\nu_{i}|(\log
(\zeta_{i}))^{\alpha}\bigg]\\&\quad\quad\quad\quad+\frac{|\mu_{2}|+|\mu_{1}|}{|\lambda|\Gamma(\alpha)}
\bigg[1+\sum_{i=1}^{n-2}|\sigma_{i}|(\log(\zeta_{i}))^{\alpha-1}\bigg]\\&\quad\quad\approx3.414437455.
\end{align*}
Now, by $(Q_{5})$ we
have,$$\quad\quad\quad\quad\quad\quad\quad\quad\quad\frac{L}{(2)((L+1)/12)(3.414437455)}
> 1 \quad\quad\quad\quad\quad\quad\quad\quad\quad\quad\quad(4.6)$$
\[\]
Hence, $L>1.320578171.$ Therefore, by Theorem3.5, the boundary value
problem (4.3) has at least one solution on $J.$

\end{document}